\documentclass{article}

\usepackage[preprint]{neurips_2025}

\usepackage[utf8]{inputenc} %
\usepackage[T1]{fontenc}    %
\usepackage{hyperref}       %
\usepackage{url}            %
\usepackage{booktabs}       %
\usepackage{amsfonts}       %
\usepackage{nicefrac}       %
\usepackage{microtype}      %
\usepackage{xcolor}         %

\usepackage{amsmath}
\usepackage{amsfonts}
\usepackage{relsize}
\usepackage{bm}
\usepackage{arydshln}
\usepackage{mathtools}
\usepackage{libertine}
\usepackage[parfill]{parskip}
\usepackage[ruled]{algorithm2e}
\SetKwComment{Comment}{/* }{ */}
\usepackage{hhline}
\usepackage{multirow}
\usepackage{subcaption}

\usepackage{fancyhdr}
\usepackage{enumitem}
\setlist[itemize]{leftmargin=*}

\let\amscases\cases
\makeatletter
\def\cases{\@ifnextchar\bgroup\plaincases\amscases}
\def\plaincases#1{\begin{cases*}#1\end{cases*}}
\makeatother

\def\eop{\hfill{$\vcenter{\hrule height1pt \hbox{\vrule width1pt height5pt
   \kern5pt \vrule width1pt} \hrule height1pt}$} \medskip}

\newtheorem{theorem}{Theorem}
\newtheorem{definition}{Definition}
\newtheorem{example}{Example}
\newtheorem{Proposition}{Proposition}
\newtheorem{Lemma}{Lemma}
\newtheorem{corollary}{Corollary}
\newtheorem{Remark}{Remark}

\newcounter{problem}

\newenvironment{proof}{\smallskip\noindent\textbf{Proof.}~}{\phantom{a}\eop}

\title{Bounding Neyman-Pearson Region \\ with $f$-Divergences}

\author{%
  Andrew Mullhaupt
     \\
  Department of Applied Mathematics and Statistics\\
  Stony Brook University\\
  Stony Brook, NY 11790 \\
  \texttt{andrew.mullhaupt@stonybrook.edu} \\
   \And
   Cheng Peng\\
  Department of Applied Mathematics and Statistics\\
  Stony Brook University\\
  Stony Brook, NY 11790 \\
   \texttt{cheng.peng.1@stonybrook.edu} \\
}

\begin{document}

\maketitle

\begin{abstract}
 
The Neyman-Pearson region of a simple binary hypothesis testing is the set of points whose coordinates represent the false positive rate and false negative rate of some  test. The lower boundary of this region is given by the Neyman-Pearson lemma, and is up to a coordinate change,  equivalent to the optimal ROC curve. 
 We establish a novel lower bound for the boundary in terms of any $f$-divergence. Since the bound  generated by hockey-stick $f$-divergences characterizes the Neyman-Pearson boundary, this bound is best possible. In the case of KL divergence, this bound improves  Pinsker's inequality.   
Furthermore, we obtain a closed-form refined upper bound for the Neyman-Pearson boundary in terms of the Chernoff 
$\alpha$-coefficient. Finally, we present methods for constructing pairs of distributions that can approximately or exactly realize any given Neyman-Pearson boundary. 
\end{abstract}

\section{Introduction}

For two simple hypotheses $H_0: X\sim P$ and $H_1: X\sim Q$, the Neyman-Pearson
region is defined by points whose coordinates represent the false positive rate and false negative rate of some hypothesis test. The Neyman-Pearson Lemma \citep{neyman1933ix} states that the likelihood ratio test achieves the smallest false positive rate  $\alpha$ for a given false negative rate  $
\beta$. This implies that the lower
boundary of the Neyman-Pearson region, termed the Neyman-Pearson boundary, is given by the likelihood ratio tests.  %

The Neyman-Pearson boundary has straightforward connections to popular classification metrics. The ROC curve of a classifier is the plot of the true positive rate against the false positive rate across all threshold levels. The Neyman-Pearson boundary is the best possible ROC curve flipped upside-down. On the other hand, we show that  the Neyman-Pearson region is  characterized by the Bayes error rate under all possible class probabilities.
Via these connections, results on ROC curve and Bayes error rate can be transferred.

Between two distributions $P$ and $Q$ one has the $f$-divergences  which
are expectations of convex functions of the likelihood ratio.  $f$-divergences provide important bounds  for the Neyman-Pearson boundary. 
Such bounds can be transparently extracted from previous literature on Bayes error rate, such as 
  \cite{kailath1967divergence,
  Chernoff1952AMO,
Hellman1970ProbabilityOE} for total variation distance, Hellinger distance, $\alpha$-divergences,  \cite{burnashev2023stein} for KL divergence, and \cite{berisha_empirically_2016-1}  for Henze-Penrose divergence.
 Our novel and general lower bound recovers the result on total variation distance, $\alpha$-divergences and KL divergence as special cases. 
Our bound generated by KL divergence adds to the literature on refinements of Pinsker's inequality
\citep{fedotov2003refinements}. Our result on tensorized bound by Chernoff $\alpha$-coefficient is closely related to the study on sample complexity of binary hypothesis testing \citep{pensia2024sample}. For a comprehensive survey on inequalities between $f$-divergences, see \cite{sason2016f}.

This study makes the following contributions: 
\begin{itemize}
\item 
 \textbf{Lower bound for Neyman-Pearson region} (Section \ref{sec_lower_bound}): 
We prove a general theorem that a function of $\alpha$ and $\beta$ is bounded by the $f$-divergence from $P$ to $Q$.  The example of the family of hockey-stick divergences leads to  the supporting lines that characterize the boundary for Neyman-Pearson region. Thus, the lower bound cannot be improved for general $f$-divergence. 
\item \textbf{Examples} (Section \ref{sec_lower_bound}): We derive several important examples of the lower bound.  The Kullback-Leibler divergence leads to a bound tighter than Pinsker's inequality.  The bound generated by $\alpha$-divergences admits tensorization, which is particularly useful in the case of multiple i.i.d. samples. 
 \item \textbf{Refined closed-form upper bound} (Section \ref{sec_upper}): We refine existing upper bounds via the convexity for Neyman-Pearson region.  We prove that the existing upper bound generated by the Chernoff $\alpha$-coefficient as well as its refined version has a closed-form expression. 
  \item  \textbf{Realization of Neyman-Pearson region} (Section \ref{sec_realization}): We show that any Neyman-Pearson boundary can be approximately realized to arbitrary precision by a pair of categorical distributions. Additionally, any Neyman-Pearson region can be exactly realized by a pair of distributions on the unit interval. 
  \item \textbf{Connection to Bayes error rate and ROC curve} (Section \ref{sec_ber}, \ref{sec_roc}): The Neyman-Pearson boundary is the best possible ROC curve flipped upside-down. The Neyman-Pearson region is  characterized by the Bayes error rate under all possible class probabilities.
Via these connections, our results directly transfer to ROC curve and Bayes error rate.

\end{itemize}

$f$-divergence has been widely used in machine learning literature for comparing distributions, such as in distributionally robust optimization \citep{ben2013robust}, generative neural network \citep{nowozin2016f},  imitation learning \citep{ke2021imitation} and AI alignment \citep{wang2023beyond}. Our results provides insight on the information contained in $f$-divergence on the Neyman-Pearson boundary, which justifies the use of $f$-divergence as constraints and objective functions from the perspective of statistical inference. 

Bayes error rate estimation has been studied in \cite{noshad_learning_2019-1, ishidaperformance, jeong2023demystifying, Theisen_Evaluating_2024} to evaluate   if models have reached optimal performance. Through our bound, estimators on various $f$-divergences \citep{nguyen2007estimating, perez2008estimation, ding2023empirical} naturally give an estimated bound on Bayes error rate. It allows us to estimate the best possible performance of a classifier before fitting a classification model.

The Wasserstein distance, which is not an $f$-divergence, is also widely
used for comparing distributions.
However, without additional restrictive assumptions, the Wasserstein
distance does not contain any 
information on statistical inference.  For a simple example, the
2-Wasserstein distance between the Gaussian distributions $\mathcal{N}(-1,\sigma)$ and
$\mathcal{N}(-1,\sigma)$ is $2$, independent of $\sigma$. But the classification problem for these
distributions is entirely determined by $\sigma$. As $\sigma$ tends to $0$, the  two
distributions get easier to tell apart, and as $\sigma$ tends to infinity, the
two distributions converge to each other in total variation distance, so
are statistically indistinguishable. So in general there is no way to
relate the classification error for these distributions to the
2-Wasserstein (or any other Wasserstein) distance. If the sample space
is bounded with respect to the metric defining a Wasserstein distance, 
then a lower bound on the total variation distance is possible \citep{gibbs2002choosing}, and  under some smoothness conditions an upper bound \citep{chae2020wasserstein}. Unlike the Wasserstein distance,  $f$-divergences are not dependent on a metric on the sample space.

\section{Preliminaries}\label{sec_prelim}

This section introduces necessary background on $f$-divergence and the Neyman-Pearson region. 
Let $P$ and $Q$ be probability measures on measurable space $(\Omega, \Sigma)$, both absolutely continuous with respect to a measure $\lambda$. Let $p(x), q(x)$ be the Radon–Nikodym derivatives, $p(x) = \frac{dP}{\mathrm{d}\lambda}$, $q(x) = \frac{dQ}{\mathrm{d}\lambda}$.
In this situation, we can use $f$-divergences to compare the distributions $P$ and $Q$.
\begin{definition}{($f$-divergence \citep{csiszar1963}.)}
Let $f: \mathsf{R}^+ \rightarrow \overline{\mathsf{R}}$ be a convex function with $f(1)=0$.  Let $f'(\infty) = \lim_{t\rightarrow 0^+} tf(\frac{1}{t})$. The  $f$-divergence from distribution $P$ to distribution $Q$ is defined by
\begin{equation}
D_f(P||Q) = \int_{q>0} f\left(\frac{p(x)}{q(x)}\right) q(x)\mathrm{d}\lambda + f'(\infty) P[q=0] .
\end{equation}
\end{definition}
The $f$-divergence does not depend on the choice of dominating measure $\lambda$. It follows from conservation of probability that for $g(x) = f(x) + k(x-1)$, $D_f(P||Q) = D_g(P||Q)$, and the converse is also true. 
Let $f^*(t) = tf\left( \frac{1}{t} \right)$. The function $f^*(t)$ is convex with $f^*(1) = 0$, and generates the divergence  $D_{f^*}(Q||P)=D_f(P||Q)$. So the $f$-divergence is symmetric if and only if $f(t) = f^*(t)$. 
Important examples of $f$ include $\frac{1}{2}|t-1|$ for total variation distance (TVD), $\frac{1}{2}(1-\sqrt{t})^2$ for squared Hellinger distance ($H^2$), and $t\log t$ for Kullback-Leibler divergence (KL).

\begin{definition}{($\alpha$-divergence \citep{havrda1967quantification, amari2000methods}.)}\label{def_alpha_div}
For $\alpha\in \mathbb{R}\backslash \{0,1\}$, the $\alpha$-divergence from $P$ to $Q$ is defined as the  $f$-divergence with $f(t) = \frac{1}{\alpha(1-\alpha)} \left( \alpha + (1-\alpha)t - t^{1-\alpha} \right)$. 
Extreme cases where $\alpha=0$ and $1$ are 
$D_{0}(P||Q) = \mathrm{KL}(P||Q)$, $D_{1}(P||Q) = \mathrm{KL}(Q||P)$. A special case is $D_{\frac{1}{2}}(P||Q) = H^2(P||Q)$. 
\end{definition}

\begin{definition}{(Chernoff $\alpha$-coefficient \citep{Chernoff1952AMO}) and connection to $\alpha$-divergence.)}\label{def_alpha_coef}
The Chernoff $\alpha$-coefficient  from $P$ to $Q$ is defined for $\alpha\in(0,1)$ by
\begin{align}
\rho_\alpha(P||Q) =& \int p(x)^\alpha q(x)^{1-\alpha} \mathrm{d} \lambda
.  \end{align}
In particular, the case of $\alpha=\frac{1}{2}$ is called the Hellinger affinity. The $\alpha$-divergence can be written as
\begin{align}
D_\alpha(P||Q) = 
\frac{1}{\alpha(1-\alpha)}\left( 1 - \rho_\alpha(P||Q) \right)
.  \end{align}
\end{definition}

To avoid confusion between the $\alpha$ in $\alpha$-divergence, Chernoff $\alpha$-coefficient and the notation for false positive rate, we use $q$ throughout the paper to represent the parameter $\alpha$ in $\alpha$-divergence and Chernoff $\alpha$-coefficient. $q \in (0,1)$.

\paragraph{Tensorization.} The Chernoff $\alpha$-coefficient for independent random variables is
\begin{equation}
\rho_\alpha \left( \bigotimes_{i=1}^n P_i ||  \bigotimes_{i=1}^n Q_i \right)
=
\prod_{i=1}^n \rho_\alpha(P_i||Q_i) .
\end{equation}
 For joint distributions of independent identically
distributed sequences this property can be written as 
$
\rho_\alpha(P^{\otimes n}||Q^{\otimes n}) = \rho_\alpha(P||Q)^n
$.
This property facilitates analysis in scenarios involving 
$n$ i.i.d. samples.

\paragraph{Binary hypothesis testing. } A hypothesis test for a pair of simple hypotheses $(H_0, H_1)$ is characterized by a function $\phi:\Omega \rightarrow [0,1]$. If $x$ is observed, $H_0$ is accepted with probability $\phi(x)$, while $H_1$ is accepted with probability $1-\phi(x)$.   The false positive rate (probability of Type I error)  $\alpha$ is the probability that $H_1$ is accepted but $H_0$ is true. The false negative rate (probability of Type II error) $\beta$ is the probability that $H_0$ is accepted but $H_1$ is true.  
A nonrandomized  hypothesis test  is characterized by a measurable set $E \in \Sigma$. If $x\in E$, $H_1$ is accepted, otherwise $H_0$ is accepted. We study the square $[0,1] \times [0,1]$ that contains all possible points $(\alpha,\beta)$.

\paragraph{Randomized test.} A randomized test randomly selects the result of test $A$   with probability $p$, and $B$ with probability $1-p$. The false positive rate and false negative rate of the randomized test are $p\alpha_a+(1-p)\alpha_b$ and $p\beta_a+(1-p)\beta_b$, where $(\alpha_a,\beta_a)$ and $(\alpha_b,\beta_b)$ are the false positive rate and false negative rate of test $A$ and $B$, respectively.

\paragraph{Neyman-Pearson region.} In the classification problem considered in this study, the hypotheses are $H_0$: $x$ is a sample from distribution $Q$, and $H_1$: $x$ is a sample from distribution $P$.  
By the Neyman-Pearson fundamental  lemma, for fixed $\alpha$, the (possibly randomized) likelihood ratio test attains the smallest $\beta$. 
For a nonrandomized test characterized by set $E$, the false positive rate $\alpha$ and false negative rate $\beta$ of the test are 
\begin{equation}
\alpha = \int_{E} dQ, \quad \beta = \int_{E^c} dP 
.  \end{equation}

\begin{definition}{(Neyman-Pearson region, Neyman-Pearson boundary, and line of ignorance.)}
The Neyman-Pearson region is the set defined by
\begin{align}
\{(\alpha,\beta): \text{there exists a  test such that the false positive rate is $\alpha$, the false negative rate is $\beta$}  \}
.  \end{align}

The line of ignorance is defined by the segment between $(0,1)$ and $(1,0)$. The tests on the line are attained by the family of randomized tests that randomly accepts the hypothesis $H_0$ with probability  $\alpha$.
\end{definition}

The Neyman-Pearson region is symmetric with respect to $(1/2, 1/2)$. The region is convex, since any point on a segment can be realized by randomizing the two tests that realizes the two ends of the segment. %
By the Neyman-Pearson fundamental  lemma, %
 the Neyman-Pearson boundary is realized by the likelihood ratio test.

\section{Lower Bound for Neyman-Pearson Region}\label{sec_lower_bound}

This section proves the novel lower bound for the Neyman-Pearson region in terms of $f$-divergence. We give examples on hockey-stick divergence, $\alpha$-divergence and KL divergence, each with important implications.

We first show  that the tests realizing extreme points of the Neyman-Pearson boundary are nonrandomized.  
A characterizing set of such tests contains almost all points that can be classified as $P$ with no error, and excludes almost all points that can be classified as $Q$ with no error. 
\begin{Lemma}%
\label{lemma_properties}
  An extreme point of the Neyman-Pearson region corresponds to a nonrandomized test. Let $E$ be the set that characterizes the nonrandomized test. We have that
\begin{equation}\label{subset}
\{x: q(x) = 0, p(x) > 0\} \in E \;\; a.e., \quad
\{x: q(x) > 0, p(x) = 0\} \in E^c \;\; a.e. 
.  \end{equation}
\end{Lemma}

Then, we prove the novel lower bound for  nonrandomized tests. 
The first step is to prove an inequality  for the false positive rate and false negative rate of nonrandomized tests in the Neyman-Pearson region. By showing that the proven inequality defines a convex body, we have that the bound is valid for randomized tests as well. 
\begin{theorem}[Convex lower bound for Neyman-Pearson region in terms of $f$-divergences.]\label{thm_inequality}
 Let $(\alpha,\beta)$ be the false positive rate and false negative rate of a nonrandomized test on the boundary for Neyman-Pearson region. Let $D_f$ be a $f$-divergence. The following inequality holds
\begin{align}\label{inequality}
(1-\alpha) f\left( \frac{\beta}{1-\alpha} \right)
+
\alpha f\left( \frac{1-\beta}{\alpha} \right)
\leq
D_f(P||Q)
.  \end{align}
The equality condition for \eqref{inequality}  is that $f(t)$ is a piecewise linear function with one singular point $t_0$, and $(\alpha,\beta)$ corresponds to the test set $E\cap \{q>0\} = \{x: p(x)/q(x) >t_0 \}$. 
The set of $(\alpha,\beta)$ satisfying \eqref{inequality}  is a convex set which includes  the Neyman-Pearson region. 
\end{theorem}

The bound is symmetric with respect to $\alpha=\beta$ if $f(t) = tf\left(\frac{1}{t}\right)$. 
The family of functions $f(t)+\lambda(t-1)$ leads to the same inequality.  
Adding a constant to $f$ leads to the same inequality. The test set corresponds to the likelihood ratio test in the Neyman-Pearson lemma. This can be viewed as a proof of the Neyman-Pearson Lemma.

We can obtain the following new bound from Theorem \ref{thm_inequality}  regarding the inverse $f$-divergence  by considering $f^*(t) = tf\left( \frac{1}{t}\right)$ as $f(t)$, or alternatively, exchanging $P$ and $Q$.
\begin{corollary}[Lower bound for Neyman-Pearson region by reversed $f$-divergences.]\label{cor_symmetric} The following bound holds
\begin{align}\label{inequality_symmetric}
\beta f\left( \frac{1-\alpha}{\beta} \right)
+
(1-\beta) f\left( \frac{\alpha}{1-\beta} \right)
\leq
D_f(Q||P)
.  \end{align}
If $f(t)=f^*(t)$, \eqref{inequality}  and \eqref{inequality_symmetric}  define the same bound , which is symmetric with respect to line $\beta= \alpha$.
\end{corollary}

The $f(t)$ that attains the equality in Theorem \ref{thm_inequality} can be written as $ \max \{ t-\gamma, 0 \}$ by  adding $\lambda(t-1)$ to $f(t)$. It belongs to the family of hockey-stick divergences, which leads to the following bound.

\begin{example}[Hockey-stick divergence.]\label{eg_hs}
Let $f(t) = \max \{ t-\gamma, 0 \}$, where $\gamma \geq 1$. The bounds obtained from \eqref{inequality} and \eqref{inequality_symmetric} can be combined to establish the following family of bounds
\begin{align}\label{eg_hs_ineq}
\beta 
\geq 
-\gamma \alpha + 1 - D_{f_\gamma}(P||Q), \quad \gamma \geq 0
.  \end{align}
\end{example}

\begin{Proposition}[Supporting lines characterizing Neyman-Pearson region.]\label{thm_supporting}
There exists a likelihood ratio test where the corresponding $(\alpha, \beta)$ makes the equality \eqref{eg_hs_ineq} hold. Varying $\gamma$ in $[0,+\infty)$, the  family of straight lines  in Example \ref{eg_hs} is the family of supporting lines to the Neyman-Pearson boundary.
\end{Proposition}
Proposition \ref{thm_supporting} shows that the bound in Theorem \ref{thm_inequality} cannot be improved for general $f$, since the family of hockey-stick divergences leads to the tightest bound. 
Two pairs distributions $P,Q$ and $P', Q'$ have the same Neyman-Pearson region if the two families of hockey-stick divergences are identical.

\begin{figure}[t]
    \centering
\begin{subfigure}[t]{0.32\textwidth}
\centering
\includegraphics[width=\textwidth]{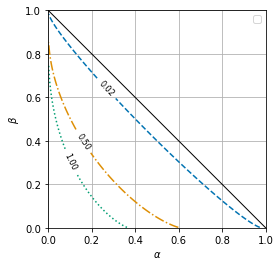}
    \caption{KL divergence}
    \label{fig_kl}
\end{subfigure}
\hfill
\begin{subfigure}[t]{0.32\textwidth}  
\centering	\includegraphics[width=\textwidth]{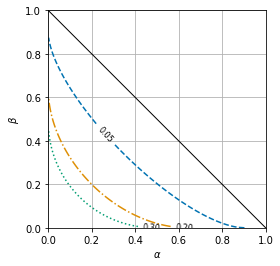}
	\caption{Squared Hellinger distance}	\label{fig_hellinger}
\end{subfigure}
\hfill
\begin{subfigure}[t]{0.32\textwidth}  
\centering
\includegraphics[width=\textwidth]{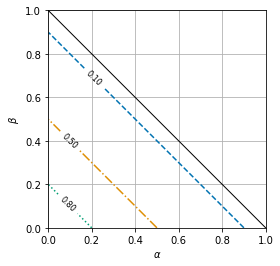}
	\caption{TVD}	\label{TVD}
\end{subfigure}
\caption{Lower bounds for Neyman-Pearson region generated by various $f$-divergences. For each $f$-divergence, $3$ different values of the divergence are used for illustration. The divergence values are labeled along the curve. The black line is the line of ignorance, which is attained by the family of randomized tests that randomly accepts the hypothesis $H_0$ with probability  $\alpha$.}\label{fig_vis}
\end{figure}

\begin{example}[Total variation distance.]\label{eg_tvd}
Let $f=\frac{1}{2}|t-1|$. %
This is a special case of hockey-stick divergence in Example \ref{eg_hs} with $\gamma = 1$. 
The bound \eqref{inequality} gives
\begin{align}
\alpha + \beta 
\geq 
1 - \mathrm{TVD}(P,Q)
.  \end{align} 
We recover the well-known bound on error rate. Substituting $f(t)$ into \eqref{inequality_symmetric}, we obtain the same bound. The equality can be attained.
\end{example}

\begin{example}[$\alpha$-divergence.]\label{eg_alpha_div}
Let $f(t) = \frac{1}{q(1-q)} \left( q + (1-q)t - t^{1-q} \right)$. The bound \eqref{inequality} gives
\begin{align}
(1-\beta)^{1-q}  \alpha^q +
\beta^{1-q}  (1-\alpha)^q
\geq
1 - q(1-q) D_q(P||Q) = \rho_q(P||Q),
\end{align}
where $\rho_q$ is the Chernoff $\alpha$-coefficient.
\end{example}

\begin{Remark}[Tensorization.]\label{remark_tensorization_1}
 To obtain the bound for independent random variables, we replace $\rho_q$ with $\prod_{i=1}^n \rho_{q, i}$ in the lower bound in Example \ref{eg_alpha_div}.  For $n$ i.i.d. samples, we simply replace $\rho$ with $\rho^n$. The tensorized version of the lower bound implies a lower bound on  the sample size  needed for a test to realize a certain $(\alpha,\beta)$.
\end{Remark}

\begin{example}[Squared Hellinger distance.]\label{eg_hellinger}
Let $f= \frac{1}{2} (\sqrt{t} -1)^2$. %
This is a special case of $\alpha$-divergence  in Example \ref{eg_hellinger} with $\alpha=\frac{1}{2}$.  
The bound \eqref{inequality} gives
\begin{align}\label{eg_hellinger_eq}
\sqrt{\alpha}\sqrt{1-\beta} + \sqrt{\beta}\sqrt{1-\alpha} \geq 1 -  H^2 =    \rho(P,Q)
.  \end{align}
where $\rho$ is the Hellinger affinity (Definition \ref{def_alpha_coef}). 
Note that $H^2(P,Q) \in [0,1]$. Substituting $f(t)$ into \eqref{inequality_symmetric}, we obtain the same bound.
\end{example}

\begin{Proposition}[Supporting lines to Hellinger lower bound.]\label{hellinger_lb_supporting}
The family of supporting lines to the bound in Example \ref{eg_hellinger} is given by
\begin{equation} \label{eq_hellinger_lb_supporting}
s\alpha + (2-s)\beta = 1 - \sqrt{1 - s(2-s)\rho^2}, \quad s \in (0, 2) .
\end{equation}
This recovers the bound in \cite{ding2023empirical}. 
\end{Proposition}

\begin{example}[Kullback-Leibler divergence.]\label{eg_KL}
Let $f = t\ln t$.  The bound \eqref{inequality} %
gives
\begin{align}\label{eg_KL_eq1}
\beta \ln \left( \frac{\beta}{1-\alpha} \right)
+
(1-\beta) \ln \left( \frac{1-\beta}{\alpha} \right)
\leq
\mathrm{KL}(P||Q)
, \end{align}
 a tighter bound than
 \begin{equation}
\alpha + \beta   \geq
 1 - 
 \sqrt{
\frac{1}{2}
\mathrm{KL}(P||Q)
}
 \end{equation}
 derived from  Pinsker's inequality
$\mathrm{TVD}(P,Q)  \leq \sqrt{\frac{1}{2}\mathrm{KL}(P||Q)}$, 
which is trivial when
$
\mathrm{KL}(P||Q)
 \geq 2$. %
\end{example}

\begin{example}[Indicator divergence.]\label{eg_ind}
Let $0 \leq \ell \leq 1 < u$. Let
\begin{align}
f(t) = \mathcal{I}_{(\ell, u)} = 
\begin{cases}
0, \quad &t \in (\ell, u) \\
+\infty, \quad &\text{otherwise}
\end{cases}
.  \end{align}
The bound obtained from \eqref{inequality} \eqref{inequality_symmetric}  is
\begin{align}\label{eg_ind_eq1}
\beta \geq \max\{-\ell \alpha + \ell, -u \alpha + 1, -u^{-1} \alpha + u^{-1}, -\ell^{-1} \alpha + 1 \}
.  \end{align}
Let $\ell = \mathrm{essinf} \frac{p(x)}{q(x)}$, $u = \mathrm{esssup}\frac{p(x)}{q(x)}$. $\max\{-\ell^{-1}, -u\}$ and $\min\{-\ell, -u^{-1}\}$ are the slopes of the tangent lines for Neyman-Pearson boundary at $(0,1)$ and $(1,0)$, respectively. 
\end{example}

\section{Upper Bound for Neyman-Pearson Boundary}\label{sec_upper}

 This section obtains the closed-form expression for the upper bound for the Neyman-Pearson boundary associated with the Chernoff $\alpha$-coefficient, and refines it using the convexity of the Neyman-Pearson region. 
 In this section, upper bounds for the Neyman-Pearson boundary refer to upper bounds for the part of the Neyman-Pearson boundary below the line of ignorance. Such an upper bound defines a convex subset of the Neyman-Pearson region.

 The Chernoff bound on BER \citep{Chernoff1952AMO,
Hellman1970ProbabilityOE} translates to  the following family of tangent lines:
 \begin{equation}\label{ber_upper_bound}
(2-s)\alpha + s\beta 
=
 s^q (2-s)^{1-q} \rho_q ,
  \quad s\in (0,2)
.  \end{equation}
 where $\rho_q$ is the Chernoff $\alpha$-coefficient (Definition \ref{def_alpha_coef}).

\begin{Proposition}[Upper bound generated by Chernoff $\alpha$-coefficients.]\label{thm_upper}
The envelope of the family of straight lines \eqref{ber_upper_bound}  is
\begin{equation}\label{upper_bound}
\beta \leq \left( q^q (1-q)^{1-q} \rho_q  \right)^{\frac{1}{q}} \alpha^{\frac{q-1}{q}},
\end{equation}
which is an upper bound of the Neyman-Pearson boundary.
\end{Proposition}

The intersection of the upper bound with the straight line $\beta = \alpha$ has coordinates $\rho_q q^q (1 - q)^{1-q}$.  
The  special case  $q=\frac{1}{2}$ results in the hyperbola $\beta \leq \frac{\rho^2}{4\alpha}$, where $\rho$ is the Hellinger affinity.

\begin{figure}[t]
\centering
\includegraphics[width=0.32\textwidth]{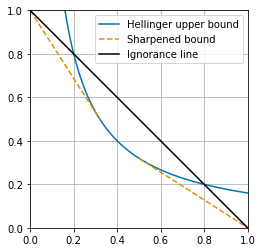}
\caption{Refined closed-form upper bound for Neyman-Pearson boundary generated by Hellinger distance. $\rho_\frac{1}{2} = 0.8$. }\label{fig_refine}
\end{figure}

Figure \ref{fig_refine} shows that upper bounds for the Neyman-Pearson boundary may lie above the line of ignorance. The line of ignorance is itself an upper bound for the Neyman-Pearson boundary. The lower convex envelope of two bounds is the weakest bound that respects them both.
\begin{theorem}[Refined upper bound.]\label{thm_refined_upper}
Any upper bound $g(\alpha)$ for the Neyman-Pearson boundary can be refined by the lower convex envelope   $\mathrm{conv}(\min\{g(\alpha), h(\alpha)\})$, where $h(\alpha)  = -\alpha + 1$ is the line of ignorance.
\end{theorem}

\begin{theorem}[Refined upper bound generated by Chernoff $\alpha$-coefficients.] \label{thm_refined_upper_closed}
Applying Theorem \ref{thm_refined_upper} to the bound \eqref{upper_bound} results in
\begin{equation}
\beta \leq \min \left\{ \left( q^q (1-q)^{1-q} \rho_q  \right)^{\frac{1}{q}} \alpha^{\frac{q-1}{q}}, \;
- \rho_q^{-\frac{1}{q}} \alpha + 1,  \;
-\rho_q^\frac{1}{q} \alpha + \rho_q^\frac{1}{q} \right\}
.  \end{equation}
\end{theorem}

\begin{Remark}[Tensorization.]
The Chernoff $\alpha$-coefficient  has a nice tensorization property for product measures (see Section \ref{sec_prelim}). The upper bound in Theorem \ref{thm_refined_upper_closed} can be generalized to $n$ i.i.d. samples by replacing $\rho_q$ with $\rho_q^n$. The tensorized upper bound indicates the lower bound on the sample size to ensure that a test can attain a desired $(\alpha,\beta)$. In practice, it can be very difficult to estimate bounds on sample size from other measure of the Neyman-Pearson region. Data-driven methods of estimating the $\alpha$-divergences make this estimation practical.
\end{Remark}

\section{Realization and Approximate Realization of the Neyman-Pearson Region}\label{sec_realization}

This section shows how pairs of distributions can be found to realize a given Neyman-Pearson region. If we pick one of the distributions to be a uniform distribution on the unit interval, the other distribution's quantile function is closely related to the Neyman-Pearson boundary. 

Consider the Neyman-Pearson region for \( P = U(0,1) \) (hence \( P \) is Lebesgue measure) and \( Q \) where the cdf of \( Q \) is \( F \). Then
\begin{equation}
\alpha([0,s]) = Q([0,s]) = F(s)
,
\quad
\beta([0,s]) = 1 - P([0,s]) = 1 - s.
\end{equation}
\begin{Lemma}\label{thm_convex_cdf}
Suppose the cdf \( F \) of a distribution \( Q \) on \( (0,1) \) is convex. Then for any measurable set \( E \) we have
\[
\alpha(E) \leq \alpha((1 - \mu(E),1)) = 1 - F(1 - \mu(E)).
\]
\end{Lemma}

\begin{theorem}\label{thm_exact_realization}
Given a Neyman-Pearson boundary $B(\alpha)$, $\alpha\in[0,1]$,  it can be realized by the pair of distributions  $P = U(0,1)$ (hence $P$ is Lebesgue measure) and $Q$ where the cdf of $Q$ is $B^{-1}(1-x)$. 
\end{theorem}

Since the Neyman-Pearson region is invariant to invertible transformation of the sample space, we can transform one of the distribution to an arbitrary distribution. For example, let $q_\mathcal{N}(p)$ be the quantile function of the standard normal distribution. The pair of distributions $q_\mathcal{N}(P)$ and $q_\mathcal{N}(Q)$ realizes the same Neyman-Pearson boundary, while $q_\mathcal{N}(Q)$ follows the standard normal distribution.

The Neyman-Pearson boundary of a pair of categorical distributions is a piecewise linear function. Naturally, we can find the piecewise linear function that is arbitrarily close to the desired Neyman-Pearson boundary, and realize it by a pair of categorical distributions.
\begin{theorem}\label{thm_approx}
Consider a  convex piecewise linear function with change points  $\{(\alpha_i, \beta_i)\}_{i=1}^n$, where $\{\alpha_i\}_{i=1}^n$ is  strictly increasing, $\{\beta_i\}_{i=1}^n$ is strictly decreasing. 
Let $(\alpha_0,\beta_0) = (0,1)$. 
Denote the negative slope by $k_i =- \frac{\beta_i - \beta_{i-1}}{\alpha_i - \alpha_{i-1}}$, $1\leq i \leq n$. If $\alpha_1 = 0$, let $k_1=+\infty$.   
Consider  a sample space of $n+1$ items and a pair of categorical distributions  represented  by $\{p_1, \cdots, p_{n+1}\}$ and $\{q_1, \cdots, q_{n+1}\}$. Without loss of generality, let $\{\frac{p_i}{q_i}\}_{i=1}^n$ be strictly decreasing. The  convex piecewise linear function is the Neyman-Pearson boundary for the pair of  categorical distributions that satisfies
\[
k_i= \frac{p_i}{q_i}, \quad i=1,\cdots, n.
\]
\end{theorem}

The two categorical distributions are specified only by the density ratios, not the densities.  We assume decreasing $\{\frac{p_i}{q_i}\}_{i=1}^n$ for notational convenience. The density ratio at each item can be arbitrarily arranged. Furthermore, the items in the sample space are arbitrary to choose. For example, they can be some positive integers.

\section{Connection to Bayes Error Rate}\label{sec_ber}

This section connects the Bayes error rate with the Neyman-Pearson boundary.  
Given the class probabilities, Bayes error rate is the lowest possible classification error under any classifiers.

\begin{definition}{(Bayes error rate.)}\label{def_ber} 
Let $(\pi_p,\pi_q)$ be the class probabilities such that $\pi_p + \pi_q = 1$. Consider a nonrandomized test characterized by set $E$. The Bayes error rate is defined by
\begin{align}
\mathrm{BER} = \min_E \left\{ \pi_q \int_E q(x) d\lambda +  \pi_p \int_{E^c} p(x) d\lambda \right\} 
.  \end{align}
\end{definition}
Denote the optimal set in Definition \ref{def_ber} by $E^*$, and the corresponding false positive rate and false negative rate by $\alpha^*$ and $\beta^*$.  Consider convex lower and upper bounds of the Neyman-Pearson boundary. %
We have $\mathrm{BER} = \pi_p \alpha^* + (1-\pi_p)\beta^*$.  

\begin{Proposition}[Bayes error rate and Neyman-Pearson region.]\label{thm_bayes} 
Consider a straight line $c = \pi_p \alpha + \pi_q \beta$.

 The following two statements are equivalent:
 
 $(i)$ $c = \pi_p \alpha + \pi_q \beta$ is a supporting line to the Neyman-Pearson boundary.
 
 $(ii)$ BER $ = c$ under class probabilities $(\pi_p, \pi_q)$.

 The following two statements are equivalent:
 
 $(iii)$  $c = \pi_p \alpha + \pi_q \beta$ is a lower bound of the Neyman-Pearson boundary.
 
 $(iv)$ BER $ \geq c$ under class probabilities $(\pi_p, \pi_q)$.

 The following two statements are equivalent:
 
 $(v)$  $c = \pi_p \alpha + \pi_q \beta$ is a supporting line to an upper bound of the Neyman-Pearson boundary.
 
 $(vi)$ BER $ \leq  c$ under class probabilities $(\pi_p, \pi_q)$.
 
\end{Proposition}

Varying $\pi_p$ in $(0,1)$, the straight line $\mathrm{BER} = \pi_p \alpha + \pi_q \beta$ represents the family of supporting lines to the Neyman-Pearson region. Thus, the study of the Neyman-Pearson region is the study of BER under all possible class probabilities.  Via this connection, classic theorems on Bayes error rate can be compared with our results.

\paragraph{Conjugate function. } 
Let $B(\alpha)$ be the minimum $\beta$ for $(\alpha,\beta)$ in the Neyman-Pearson region. %
 Let $B^*(z)$ be the convex conjugate of $B(\alpha)$
\begin{equation}\label{conjugate}
B^*(z) = \min_{\alpha\in[0,1]} \{ z\alpha - B(\alpha)  \}.
\end{equation}
Given $z$, the minimum is attained at $\alpha'$ such that $z \in \partial B(\alpha')$. We have $B^*(z) = z \alpha' - B(\alpha')$. Since $(\alpha', B(\alpha'))$ is on $B(\alpha)$ and $z \in \partial B(\alpha')$, we have the supporting line at $\alpha'$, $B^*(z) = z\alpha - \beta$. From Proposition \ref{thm_bayes},  $\pi_p \alpha + \pi_q \beta = \mathrm{BER}$ is the family of supporting lines. Solving for $\pi_p$ and BER,
\begin{equation}
\pi_p = \frac{z}{z-1}, \quad
\mathrm{BER} = (\pi_p - 1)	B^*(z).
\end{equation}
Thus we obtain the connection between BER and the slope of the supporting line of the Neyman-Pearson region
\begin{equation}
\mathrm{BER} = \frac{B^*(z)}{z-1}.
\end{equation}

\section{Connection to ROC Curve}\label{sec_roc}

This section discusses the connection between the Neyman-Pearson boundary and the ROC curve in binary classification, and proves that any ROC curve can be realized by a one-parameter family  of tests that randomizes the tests on the Neyman-Pearson boundary.

\paragraph{Equivalence to best possible ROC curve.} The ROC curve of a classifier is the plot of the true positive rate  against the false positive rate across all threshold levels. Since the true and false positive rates sum to $1$, the ROC curve has a one-to-one correspondence to a curve of false positive rate and false negative rate considered in this study.   When the classifier is applied to the population, the corresponding  curve of false positive rate and false negative rate must be  contained in the Neyman-Pearson region according to the Neyman-Pearson fundamental  lemma. Thus, our result can be interpreted as bounding the best possible ROC curve in binary classification with $f$-divergences.

Next, we realize an arbitrary ROC curve by randomizing tests on the Neyman-Pearson boundary. Excluding a possibly vertical part at  $\alpha=0$, the Neyman-Pearson boundary on interval $(0,1]$ is a convex function. We denote it by $f(\alpha)$. An ROC curve is arbitrary in shape. It can be noncontinuous or nonconvex, as long as it is a function defined on $(0,1]$.

\begin{Proposition}[Realization of ROC Curve by One-Parameter Family of Tests]\label{thm_roc_realization}
Consider an arbitrary ROC curve $g(t)$ parameterized by one parameter $t$. 
Given the tests that realizes the Neyman-Pearson boundary, the  ROC curve  can be realized by a one-parameter family of tests that randomizes the tests on the Neyman-Pearson boundary and the tests on the line of ignorance.
\end{Proposition}

\begin{figure}[t]
\centering
\begin{minipage}[t]{0.32\textwidth}
  \centering
  \includegraphics[width=\linewidth]{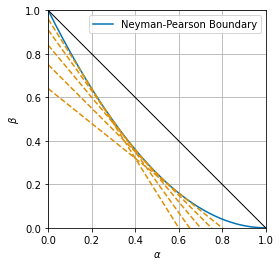}
  \caption{Supporting lines to Neyman-Pearson region corresponding to hockey-stick divergences with different $\gamma$.}
  \label{fig_hs}
\end{minipage}\hfill
\begin{minipage}[t]{0.32\textwidth}
  \centering
  \includegraphics[width=\linewidth]{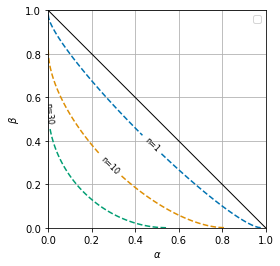}
  \caption{Lower bounds for Neyman-Pearson region corresponding to squared Hellinger distance with $n$ samples. $\rho = 0.99$.}
  \label{fig_Hellinger_n}
\end{minipage}\hfill
\begin{minipage}[t]{0.32\textwidth}
  \centering
  \includegraphics[width=\linewidth]{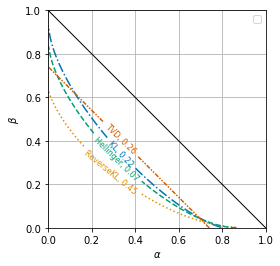}
  \caption{Lower bounds for Neyman-Pearson region for various divergences between concentric univariate Gaussians.}
  \label{fig_div_pair}
\end{minipage}
\end{figure}

\section{Visualization}\label{sec_vis}

\paragraph{Shape of lower bounds.}
Figure \ref{fig_vis} shows the shape of the lower bounds corresponding to  TVD (Example \ref{eg_tvd}),  Hellinger divergence (Example \ref{eg_hellinger}) and KL divergence (Example \ref{eg_KL}). For each divergence, we visualize the bound using different divergence values. The larger the divergence, the easier it is to distinguish the two distributions, the farther the lower bound is from the line of ignorance. The lower bound generated by KL-divergence is not symmetric.

\paragraph{Supporting lines via hockey-stick divergences.}
Figure \ref{fig_hs} shows several supporting lines in the family obtained in Example \ref{eg_hs}. The family of supporting lines characterizes the Neyman-Pearson boundary. The Neyman-Pearson boundary can be generated by a uniform distribution and a beta distribution $\mathrm{Beta}(1, 1/2)$ on the unit interval. As we show in Section \ref{sec_realization}, the boundary is not uniquely associated with the pair of distributions.

\paragraph{Tensorized Hellinger bound for product measure.}
Figure \ref{fig_Hellinger_n} shows the lower bounds corresponding to the Hellinger distance with different sample size (see Remark \ref{remark_tensorization_1}).  Lower bound of larger sample size are farther from the line of ignorance. The larger the sample size, the easier it is to tell which distribution the samples are drawn from.

\paragraph{Shape of bounds for concentric Gaussian distribution pair.}
Figure \ref{fig_div_pair} shows the lower bounds corresponding to various divergences for a given pair of concentric Gaussian distributions. 
 The bounds differ in shape and intersect one another, showing that no bound uniformly dominates all others across $\alpha$.
 
\paragraph{Refined Hellinger bound.}
Figure \ref{fig_refine} shows the closed-form upper bound in Proposition \ref{thm_upper} and the refined upper bound in Theorem \ref{thm_refined_upper_closed}, where $q = \frac{1}{2}$. 
 The middle part of the original bound is retained, while the rest of the  bound are refined by two segments.

\section{Concluding Remarks}

We advance the   understanding of the Neyman-Pearson region by establishing a general lower bound in terms of  $f$-divergences and deriving important examples. We also provide refined upper bounds, including a closed-form expression based on the Chernoff $\alpha$-coefficient. Furthermore, the realizability results and the connections to Bayes error rate and ROC curves sets the stage for broader applications in statistical decision theory.


\bibliography{NP_ref.bib}
\bibliographystyle{apalike}

\appendix

\section{Proofs}

\subsection{Proof of Lemma \ref{lemma_properties}}

\begin{proof}
First, we prove that the test is nonrandomized. The Neyman-Pearson region is a convex set in $[0,1]\times [0,1]$ that contains $(0,1)$ and $(1,0)$. A randomized test can be represented by a point in the segment between two points corresponding to two tests. Since an extreme point is not in the middle of any segment in the convex set, a test corresponding to an extreme point is nonrandomized.

Next, we prove \eqref{subset}. Let $(\alpha,\beta)$ be the false positive rate and false negative rate of the test characterized by $E$. Suppose $\lambda(\{x: q(x) = 0, p(x) > 0\} \backslash  E) > 0$. Consider a new test characterized by $E \cup \{x: q(x) = 0, p(x) > 0\}$ with false positive rate and false negative rate $(\hat{\alpha},\hat{\beta})$. The new test has a  strictly smaller false negative rate, since $p(x)>0$. The new test has the same false positive rate since $q(x) = 0$. $(\alpha,\beta)$ is in the interior of any convex set that contains $(0,1), (1,0), (\hat{\alpha},\hat{\beta})$.  This contradicts the condition that $(\alpha,\beta)$ is an extreme point of the Neyman-Pearson region. Thus, we proved $\{x: q(x) = 0, p(x) > 0\} \in E \;\; a.e.$. The second part of  \eqref{subset} can be proved similarly.
\end{proof}

\subsection{Proof of Theorem \ref{thm_inequality}}\label{proof_thm_ineq}

\begin{proof}
Let $E$ be the set that characterizes the considered nonrandomized test. 
We prove two inequalities, whose sum gives the desired result
\begin{align}
(1-\alpha) f\left( \frac{\beta}{1-\alpha} \right)
&\leq
\int_{E^c  \cap \{q>0\} } f\left( \frac{p(x)}{q(x)} \right)  q(x)  \mathrm{d}\lambda
+
f'(\infty)P[E^c \cap  \{q=0\}]\label{lb_eq1}
\\
\alpha f\left( \frac{1-\beta}{\alpha} \right)
&\leq
\int_{E \cap  \{q>0\}} f\left( \frac{p(x)}{q(x)} \right)  q(x)  \mathrm{d}\lambda
+
f'(\infty)P[E \cap  \{q=0\}]\label{lb_eq2}
.  \end{align}
For the cases of $\alpha=0$ and $\alpha=1$, we adopt the convention that $0\cdot f(\cdot) = 0$.

First, we prove \eqref{lb_eq1}. 
By Lemma \ref{lemma_properties}, $P[E^c \cap  \{q=0\}] = 0$. We have
\begin{align}
\int_{E^c \cap \{q>0\}} f\left( \frac{p(x)}{q(x)} \right)  q(x) \mathrm{d}\lambda 
&= 
(1-\alpha)\int_{E^c  \cap \{q>0\}} f\left( \frac{p(x)}{q(x)} \right)  \frac{1}{1-\alpha}  q(x)  \mathrm{d}\lambda \\
& \geq
(1-\alpha)f \left( \int_{E^c \cap \{q>0\}}   \frac{p(x)}{q(x)}  \frac{1}{1-\alpha} q(x) \mathrm{d}\lambda \right) \label{Jensen}
\\
& =
(1-\alpha)f \left( \int_{E^c}  p(x)   \frac{1}{1-\alpha} \mathrm{d}\lambda \right)\label{inequality_pre}
\\
&=
(1-\alpha)f\left( \frac{\beta}{1-\alpha}  \right) \label{inequality_part}
.  \end{align}
Note that since $\frac{p(x)}{q(x)}$ is measurable on $E^c \cap \{q > 0\}$, we can apply Jensen's inequality  to obtain \eqref{Jensen}. %

To prove \eqref{lb_eq2}, we use \eqref{lb_eq1}. 
Consider $f^*(t) = tf\left(\frac{1}{t}\right)$. $f^*(t)$ is also a convex function. Consider a new pair of hypotheses  $\tilde{H}_0:$ $x$ is a sample from $Q$, and $\tilde{H}_1:$ $x$ is a sample from $P$. Consider the test characterized by $E^c$. Notice that $\tilde{\alpha} = 1-\alpha$, $\tilde{\beta} = 1-\beta$. By \eqref{lb_eq1}, we have
\begin{align}
\int_{E\cap \{p>0\}} f^*\left( \frac{q(x)}{p(x)} \right)  p(x) \mathrm{d}\lambda 
+
f^{*'}(\infty) Q[E \cap \{p=0\}]
&\geq 
(1-\tilde{\alpha})f^*\left( \frac{\tilde{\beta}}{1-\tilde{\alpha}}  \right) 
.  \end{align}
By Lemma \ref{lemma_properties},  $Q[E \cap \{p=0\}]=0$. Thus
\begin{align}
\int_{E\cap \{p>0\}} f\left( \frac{p(x)}{q(x)} \right)  q(x) \mathrm{d}\lambda 
&\geq 
\alpha f\left( \frac{1-\beta}{\alpha}  \right) 
.  \end{align}

With \eqref{lb_eq1}\eqref{lb_eq2} proved, we have proved \eqref{inequality}.

Since $f(t)$ is convex, its perspective function $t f\left( \frac{s}{t}\right)$ is convex. Thus the left-hand side of \eqref{inequality} is convex in $(\alpha,\beta)$. Since the level set of a convex function is convex, \eqref{inequality} defines a convex set.

Since a randomized test can be represented by a point in the segment connecting two points that correspond to two tests, the convex hull of points of nonrandomized test contains all points of randomized test. Since the convex hull is a subset of any convex set that contains the points of nonrandomized test, the convex set defined by \eqref{inequality} contains the Neyman-Pearson region. The boundary below the line of ignorance $\beta =-\alpha+1$ is a lower bound of the Neyman-Pearson region.
\end{proof}

\subsection{Proof of Corrollary \ref{cor_symmetric}}

\begin{proof}
Consider $f^*(t) = tf\left( \frac{1}{t}\right)$ as $f(t)$ in \eqref{inequality},
\begin{align}
(1-\alpha) f^*\left( \frac{\beta}{1-\alpha} \right)
+
\alpha f^*\left( \frac{1-\beta}{\alpha} \right)
\leq
D_{f^*}(P||Q)
.  \end{align}
Thus
\begin{align}
\beta f\left( \frac{1-\alpha}{\beta} \right)
+
(1-\beta) f\left( \frac{\alpha}{1-\beta} \right)
\leq
D_f(Q||P)
.  \end{align}
\end{proof}

\begin{Remark}
An alternative proof of Corollary \ref{cor_symmetric} is provided below.

\begin{proof}
Switch $P$ and $Q$ in \eqref{inequality},
\begin{align}
(1-\tilde{\alpha}) f\left( \frac{\tilde{\beta}}{1-\tilde{\alpha}} \right)
+
\tilde{\alpha} f\left( \frac{1-\tilde{{\beta}}}{\tilde{\alpha}} \right)
\leq
D_f(Q||P)
.  \end{align}
Since $P$ and $Q$ are switched, $\alpha = \tilde{\beta}$, $\beta = \tilde{\alpha}$. 
\end{proof}
\end{Remark}

\subsection{Proof of Example \ref{eg_hs}}

\begin{proof} 
Substituting $f(t)$ into \eqref{inequality}, 
\begin{align}
(1-\alpha) \max\left\{  \frac{\beta}{1-\alpha} - \gamma, 0 \right\}
+
\alpha \max\left\{  \frac{1-\beta}{\alpha} - \gamma, 0 \right\}
\leq
D_f(P||Q)
= 
\int  \max\left\{\frac{p}{q}-\gamma, 0\right\}q \mathrm{d}\lambda
.  \end{align}
\begin{align}
 \max\left\{  \beta + \gamma \alpha - \gamma , 0 \right\}
+
\max\left\{  1-\beta- \gamma \alpha , 0 \right\}
\leq
\int  \max\left\{p-\gamma q, 0\right\} \mathrm{d}\lambda
=
\int_{ p-\gamma q > 0}  (p-\gamma q) \mathrm{d}\lambda
.  \end{align}
Since $\beta + \gamma \alpha < \gamma$ for the points in the triangle with vertices $(0,0), (0,1), (1,0)$ in consideration, the first term disappears. We have
\begin{align}\label{proof_hs_1}
\beta 
\geq 
-\gamma \alpha + 1 - D_{f_\gamma}(P||Q).
.  \end{align}

Substituting $f(t)$ into \eqref{inequality_symmetric}, we have
\begin{align}
\max\{1-\alpha-\gamma\beta , 0\} +
\max\{\alpha - \gamma +\gamma  \beta, 0 \}
 \leq 
 \int \max\{q - \gamma p, 0\} \mathrm{d}\lambda
 =
  \int_{q - \gamma p > 0} (q - \gamma p) \mathrm{d}\lambda
.  \end{align}
Since $\beta < - \frac{1}{\gamma} \alpha + 1$  for the points in the triangle with vertices $(0,0), (0,1), (1,0)$ in consideration, the second term disappears. We have
\begin{align}
\beta \geq -\frac{1}{\gamma}\alpha + \frac{1}{\gamma}(1-D_{f_\gamma}(Q||P)) 
.  \end{align}

To prove that the two families of straight lines can be combined and get \eqref{eg_hs_ineq}, we need the following relation
\begin{align}
D_{f_\gamma}(Q||P)  =& \int  \max\left\{q-\gamma p, 0\right\} \mathrm{d}\lambda
=
\int_{ q-\gamma p > 0}  (q-\gamma p) \mathrm{d}\lambda 
\\
D_{f_{\frac{1}{\gamma}}}(P||Q)  =&
\int  \max\left\{p- \frac{1}{\gamma} q, 0\right\} \mathrm{d}\lambda
=
\int_{ p- \frac{1}{\gamma} q > 0}  (p-\frac{1}{\gamma} q) \mathrm{d}\lambda
=
 \frac{1}{\gamma} \int_{ \gamma p -  q > 0}  (\gamma p- q) \mathrm{d}\lambda
 \\
=& 
 \frac{1}{\gamma} \left( \int (\gamma p - q) d\lambda -  \int_{ \gamma p -  q < 0}  (\gamma p- q) \mathrm{d}\lambda \right)
 \\
=&
 1 - \frac{1}{\gamma} + \frac{1}{\gamma}  D_{f_\gamma}(Q||P)
.  \end{align}
Replace $\gamma$ with $\frac{1}{\gamma}$ in \eqref{proof_hs_1}, we have
\begin{align}
\beta
 \geq&
 - \frac{1}{\gamma} \alpha + 1 - D_{f_\frac{1}{\gamma}}(P||Q)
 \\
=&
 - \frac{1}{\gamma} \alpha + 1 - 1 + \frac{1}{\gamma} -  \frac{1}{\gamma} D_{f_\gamma}(Q||P)
 \\
 =&
 -\frac{1}{\gamma}\alpha + \frac{1}{\gamma}(1-D_{f_\gamma}(Q||P)) 
.  \end{align}
Note that $\frac{1}{\gamma} \in [0,1]$. This shows that the two families are included in one by letting $\gamma$ in \eqref{eg_hs_ineq} be in $[0,+\infty)$.
\end{proof}

\subsection{Proof of Proposition \ref{thm_supporting}}

\begin{proof}
In the proof of Theorem \ref{thm_inequality}, the equality in Jensen's inequality holds if and only if $f(t)$ is linear on the set $E\cap\{q>0\}$. In Example  \ref{eg_hs}, $f(t)$ is a piecewise linear function with breakpoint $t=\gamma$. 
 Thus, the equality in \eqref{eg_hs_ineq} holds if and only if the following two conditions hold: $(i)$ $E\cap\{q>0\} \subset \{x: p(x)/q(x) > \gamma\}$ or $E\cap\{q>0\} \subset \{x: p(x)/q(x) < \gamma\}$; $(ii)$ $E^c\cap\{q>0\} \subset \{x: p(x)/q(x) > \gamma\}$ or $E^c\cap\{q>0\} \subset \{x: p(x)/q(x) < \gamma\}$. 
Then, the equality in \eqref{eg_hs_ineq} holds if and only if $E\cap\{q>0\} = \{x: p(x)/q(x) > \gamma\}$ or $E\cap\{q>0\} = \{x: p(x)/q(x) < \gamma\}$. Each set corresponds to a point in the Neyman-Pearson region. The two points are symmetric with respect to $(\frac{1}{2}, \frac{1}{2})$. We are concerned with the first point, which is in the triangle with vertices $(0,0), (0,1), (1,0)$. The set corresponds to the likelihood ratio test with threshold $\gamma$. 

Since we considered the equality condition, the point $(\alpha,\beta)$  is the only intersection between the straight line and the Neyman-Pearson region.  Thus every straight line in the family is a supporting line to the Neyman-Pearson region.

Varying $\gamma$ in $[0,+\infty)$, we obtain a family of supporting lines. Since the family contains lines with any non-positive slope, every supporting line to the Neyman-Pearson region must belong to this family.
\end{proof}

\subsection{Proof of Example \ref{eg_tvd}}

\begin{proof}
\begin{align}
\frac{1}{2}|\beta - 1 + \alpha|
+ \frac{1}{2}|1-\beta - \alpha| \leq \mathrm{TVD}(P,Q)
.  \end{align}
\end{proof}

\subsection{Proof of Example \ref{eg_hellinger}}

\begin{proof}
\begin{align}
\frac{1}{2}  (1-\alpha)\left( \sqrt{\frac{\beta}{1-\alpha}}-1 \right)^2
+
\frac{1}{2}  \alpha \left( \sqrt{\frac{1-\beta}{\alpha}}-1 \right)^2
\leq H^2(P,Q).
.  \end{align}
\begin{align}
\beta + 1 - \alpha - 2\sqrt{\beta(1-\alpha)}
+
1-\beta + \alpha - 2\sqrt{(1-\beta)\alpha}
\leq 2 H^2(P,Q).
.  \end{align}
\end{proof}

\begin{Remark}
Example 2 has the following alternative proof.

To prove
\begin{align}
\left( \sqrt{\int_{E^c} p(x) \mathrm{d}\lambda} - \sqrt{\int_{E^c} q(x)\mathrm{d}\lambda} \right)^2
+
\left( \sqrt{\int_{E} p(x)\mathrm{d}\lambda} - \sqrt{\int_{E} q(x)\mathrm{d}\lambda} \right)^2
\leq
\int (\sqrt{q(x)}-\sqrt{p(x)})^2 \mathrm{d}\lambda
.  \end{align}
we need to prove
\begin{align}
2 - 2\sqrt{\int_{E} p(x) \mathrm{d}\lambda \int_{E} q(x)\mathrm{d}\lambda} - 2\sqrt{\int_{E^c} p(x) \int_{E^c} q(x)\mathrm{d}\lambda}
\leq
2 - 2\int \sqrt{p(x)q(x)} \mathrm{d}\lambda
.  \end{align}
which is to prove
\begin{align}
\sqrt{\int_{E} p(x) \mathrm{d}\lambda\int_{E} q(x)\mathrm{d}\lambda } + \sqrt{\int_{E^c} p(x) \mathrm{d}\lambda \int_{E^c} q(x) \mathrm{d}\lambda}
\geq
\int_{E} \sqrt{p(x)q(x)} \mathrm{d}\lambda + \int_{E^c} \sqrt{p(x)q(x)} \mathrm{d}\lambda
.  \end{align}
This is true by Cauchy Schwartz inequality.
\end{Remark}

\subsection{Proof of Proposition \ref{hellinger_lb_supporting}}

\begin{proof}
First, we simplify the curve  using a trigonometric substitution to to a simpler form. We then derive the general equation for the tangent line to the curve based on this simplified form, using a parameter $\theta$. Finally, we show that this derived tangent line equation is equivalent to the given family of lines, by establishing a direct relationship between the curve parameter $\theta$ and the line family parameter $s$.

Let $\alpha = \sin^2 \theta$ and $\beta = \sin^2 \phi$, where $\theta, \phi \in (0, \pi/2)$. Substituting into the bound \eqref{eg_hellinger_eq}
\begin{align}
\rho =& \sqrt{\sin^2 \theta (1-\sin^2 \phi)} + \sqrt{\sin^2 \phi (1-\sin^2 \theta)}  
= \sin(\theta + \phi) . \label{eq_trig1}
\end{align}
Let $\gamma = \arcsin \rho$,  $\gamma \in (0, \pi/2)$. \eqref{eq_trig1} has two solutions for $\theta + \phi \in (0, \pi)$, $\theta + \phi = \gamma$ and $\theta + \phi = \pi - \gamma$. 
We consider the first case. We have $\theta \in (0, \gamma)$.

The curve \eqref{eg_hellinger_eq} can be parameterized by $\theta$
\begin{align}
\alpha(\theta)= \sin^2 \theta , \quad
\beta(\theta) = \sin^2(\gamma - \theta).
\end{align}
Take derivative with respect to $\theta$ 
\begin{align}
\frac{d\alpha}{d\theta} = 2 \sin \theta \cos \theta = \sin(2\theta) , \quad
\frac{d\beta}{d\theta} = -\sin(2(\gamma - \theta))
\end{align}
The slope of the tangent line is 
\begin{align}
m =  \frac{d\beta/d\theta}{d\alpha/d\theta} = \frac{-\sin(2\gamma - 2\theta)}{\sin(2\theta)}
\end{align}
The tangent line equation at $(\alpha_0, \beta_0) = (\sin^2 \theta, \sin^2(\gamma - \theta))$ is
\begin{align}
\beta - \sin^2(\gamma - \theta) &= \frac{-\sin(2\gamma - 2\theta)}{\sin(2\theta)} (\alpha - \sin^2 \theta) \\
\beta \sin(2\theta) - \sin(2\theta) \sin^2(\gamma - \theta) &= -\alpha \sin(2\gamma - 2\theta) + \sin(2\gamma - 2\theta) \sin^2 \theta \\
\alpha \sin(2\gamma - 2\theta) + \beta \sin(2\theta) &= \sin(2\gamma - 2\theta) \sin^2 \theta + \sin(2\theta) \sin^2(\gamma - \theta)
\end{align}
Simplify the right-hand side (RHS)
\begin{align}
\text{RHS} &= (2 \sin(\gamma - \theta) \cos(\gamma - \theta)) \sin^2 \theta + (2 \sin \theta \cos \theta) \sin^2(\gamma - \theta) \\
&= 2 \sin \theta \sin(\gamma - \theta) (\cos(\gamma - \theta) \sin \theta + \cos \theta \sin(\gamma - \theta)) \\
&= 2 \sin \theta \sin(\gamma - \theta) (\sin(\theta + (\gamma - \theta))) \\
&= 2 \sin \theta \sin(\gamma - \theta) \sin \gamma \\
&= 2 \rho \sin \theta \sin(\gamma - \theta).
\end{align}
Thus, the tangent line equation is
\begin{equation} \label{eq_tangent}
\alpha \sin(2\gamma - 2\theta) + \beta \sin(2\theta) = 2 \rho \sin \theta \sin(\gamma - \theta).
\end{equation}

Next, we compare \eqref{eq_tangent} with \eqref{eq_hellinger_lb_supporting}. 
We seek a scaling factor $k$ and a parameter $s(\theta)$ such that
\begin{align}
s &= k \sin(2\gamma - 2\theta) \label{eq_s_k} \\
2-s & = k \sin(2\theta) \label{eq_2ms_k} \\
1 - \sqrt{1 - s(2-s)\rho^2} & = k (2 \rho \sin \theta \sin(\gamma - \theta)) \label{eq_C_k}
\end{align}
Adding \eqref{eq_s_k} and \eqref{eq_2ms_k}
\begin{align}
2 =& k (\sin(2\gamma - 2\theta) + \sin(2\theta)) \\
=& 2k \sin(\tfrac{2\gamma - 2\theta + 2\theta}{2}) \cos(\tfrac{2\gamma - 2\theta - 2\theta}{2}) \\
=& 2k \sin \gamma \cos(\gamma - 2\theta) \\
=& 2k \rho \cos(\gamma - 2\theta) 
\end{align}
Thus, 
\begin{align}
k &= \frac{1}{\rho \cos(\gamma - 2\theta)}.
\end{align}
Note that for $\theta \in (0, \gamma)$ and $\gamma \in (0, \pi/2)$, $\gamma - 2\theta \in (-\gamma, \gamma) \subset (-\pi/2, \pi/2)$, so $\cos(\gamma - 2\theta) > 0$.
Substituting $k$ back into \eqref{eq_s_k}
\begin{align}
s = s(\theta) = \frac{\sin(2\gamma - 2\theta)}{\rho \cos(\gamma - 2\theta)}
\end{align}

Lastly, we need to verify  \eqref{eq_C_k}. The right-hand side of \eqref{eq_C_k} is
\begin{align}
\frac{1}{\rho \cos(\gamma - 2\theta)} (2 \rho \sin \theta \sin(\gamma - \theta)) = \frac{2 \sin \theta \sin(\gamma - \theta)}{\cos(\gamma - 2\theta)}.
\end{align}
We study the left-hand side in multiple steps for clarity. 

First, we  evaluate the term under the square root
\begin{align}
s(2-s)\rho^2 &= (k \sin(2\gamma - 2\theta)) (k \sin(2\theta)) \rho^2 \\
&= k^2 \rho^2 \sin(2\gamma - 2\theta) \sin(2\theta) \\
&= \left(\frac{1}{\rho \cos(\gamma - 2\theta)}\right)^2 \rho^2 \sin(2\gamma - 2\theta) \sin(2\theta) \\
&= \frac{\sin(2\gamma - 2\theta) \sin(2\theta)}{\cos^2(\gamma - 2\theta)}
\end{align}
Then,
\begin{align}
1 - s(2-s)\rho^2
 = 1 - \frac{\sin(2\gamma - 2\theta) \sin(2\theta)}{\cos^2(\gamma - 2\theta)}
 = \frac{\cos^2(\gamma - 2\theta) - \sin(2\gamma - 2\theta) \sin(2\theta)}{\cos^2(\gamma - 2\theta)}.
\end{align}
Given that
$
\sin(2\gamma - 2\theta) \sin(2\theta) 
= \frac{1}{2}\big( \cos(2\gamma - 4\theta) - \cos(2\gamma)\big) $, 
the numerator is
\begin{align*}
&\cos^2(\gamma - 2\theta) - \sin(2\gamma - 2\theta) \sin(2\theta)\\
 =& \cos^2(\gamma - 2\theta) - \frac{1}{2}[\cos(2(\gamma - 2\theta)) - \cos(2\gamma)] \\
=& \cos^2(\gamma - 2\theta) - \frac{1}{2}[ (2\cos^2(\gamma - 2\theta) - 1) - \cos(2\gamma)]  \\
=& \cos^2(\gamma - 2\theta) - \cos^2(\gamma - 2\theta) + \frac{1}{2} + \frac{1}{2}\cos(2\gamma) \\
=&  \cos^2\gamma 
\end{align*}
Thus, 
\begin{align}
1 - \sqrt{1 - s(2-s)\rho^2}
 = 
 1 - \frac{\cos\gamma}{\cos(\gamma - 2\theta)}
\end{align}
Putting the left-hand side and right-hand side of \eqref{eq_C_k} together, 
we need to check
\begin{align}
1 - \frac{\cos\gamma}{\cos(\gamma - 2\theta)} = \frac{2 \sin \theta \sin(\gamma - \theta)}{\cos(\gamma - 2\theta)},
\end{align}
which is to check
\begin{align}
\cos(\gamma - 2\theta) - \cos\gamma = 2 \sin \theta \sin(\gamma - \theta).
\end{align}
This can be  verified by the product-to-sum equality. The proof is completed.

We have shown that for every parameter $s \in (0, 2)$, the line $s\alpha + (2-s)\beta = 1 - \sqrt{1 - s(2-s)\rho^2}$ corresponds to a supporting  line of the curve $\sqrt{\alpha(1-\beta)} + \sqrt{\beta(1-\alpha)} = \rho$ at the point $(\alpha(\theta), \beta(\theta)) = (\sin^2\theta, \sin^2(\gamma-\theta))$ where $\theta$ is uniquely determined by $s = \frac{\sin(2\gamma - 2\theta)}{\rho \cos(\gamma - 2\theta)}$.
\end{proof}

\subsection{Proof of Example \ref{eg_ind}}

\begin{proof}
If $\ell < \mathrm{essinf} \frac{p(x)}{q(x)}$ and $u > \mathrm{esssup} \frac{p(x)}{q(x)}$, then 
 $D_f(P||Q) = 0$. Finiteness of $f\left( \frac{\beta}{1-\alpha} \right)$ and $f\left( \frac{1-\beta}{\alpha} \right)$ requires
\begin{align}
\ell &\leq \frac{\beta}{1-\alpha} \leq u \\
\ell &\leq \frac{1-\beta}{\alpha} \leq u 
.  \end{align}
That is,
\begin{align}
\beta &\leq -u \alpha + u \label{proof_indi_div_1} \\
\beta &\geq -\ell \alpha + \ell \label{proof_indi_div_2} \\
\beta &\leq -\ell \alpha + 1 \label{proof_indi_div_3} \\
\beta &\geq -u \alpha + 1  \label{proof_indi_div_4}
.  \end{align}
Since $0 \leq \ell \leq 1 < u$, \eqref{proof_indi_div_1} and  \eqref{proof_indi_div_3} hold for all points in the considered triangle with vertices $(0,0), (0,1), (1,0)$. \eqref{proof_indi_div_2} and \eqref{proof_indi_div_4} lead to \eqref{eg_ind_eq1}. \eqref{eg_ind_eq2} can be proved similarly.
\end{proof}

\subsection{Proof of Example \ref{eg_KL}}
\begin{proof}
Let $s = \sqrt{\frac{1}{2}\mathrm{KL}}$.  We prove that the straight line $\beta = -\alpha  + 1 - s$ has at most one intersection with the convex curve \eqref{eg_KL_eq1}. Define the following $G(\alpha,s)$ by replacing $\beta$ in \eqref{eg_KL_eq1} with $-\alpha  + 1 - s$ and moving $\mathrm{KL}(P||Q)$ to the same side. We prove that $G(\alpha, s)$ has no more than one root.  
\begin{align}
G(\alpha, s) = 
(-\alpha + 1-s) \ln \left( \frac{-\alpha + 1-s}{1-\alpha} \right)
+
(\alpha + s) \ln \left( \frac{\alpha + s}{\alpha} \right)
 - 2s^2
.  \end{align}
\begin{align}
\frac{\partial G}{\partial s}
=
-4 s + \log\frac{\alpha + s}{\alpha} - \log \frac{-1 + \alpha +s}{-1 + \alpha}
.  \end{align}
\begin{align}
\frac{\partial^2 G}{\partial s^2}
=
-4 + \frac{1}{1 - \alpha - s} + \frac{1}{\alpha + s} \geq 0
.  \end{align}
The equailty holds iff $\alpha + s = \frac{1}{2}$. Thus, $G(\alpha, s)$ is convex in $s$. 
Since $\frac{\partial G}{\partial s}(\alpha,0) = 0$,  $G(\alpha, s)$ is increasing in $s$ for $s>0$. Since $G(\alpha, 0) = 0$,  we have $G\geq 0$ for $s>0$.
\end{proof}

\subsection{Proof of Example \ref{eg_ind}}

\begin{proof}
If $\ell > \mathrm{essinf} \frac{p(x)}{q(x)}$ or $u < \mathrm{esssup} \frac{p(x)}{q(x)}$ , $D_f(P||Q) = +\infty$. The bound holds trivially.

If $\ell < \mathrm{essinf} \frac{p(x)}{q(x)}$ and $u > \mathrm{esssup} \frac{p(x)}{q(x)}$, substituting $f(t)$ into \eqref{inequality}, we have
\begin{align}\label{eg_ind_eq1}
\beta \geq \max\{-\ell \alpha + \ell, -u \alpha + 1 \}
.  \end{align}
Substituting $f(t)$ into \eqref{inequality_symmetric}, we have
\begin{align}\label{eg_ind_eq2}
\beta \geq \max\{-u^{-1} \alpha + u^{-1}, -\ell^{-1} \alpha + 1 \}
.  \end{align}
\end{proof}

\subsection{Example of $\chi^2$-divergence}

\begin{example}[$\chi^2$-divergence.] Let $f(t) = (t-1)^2$. Substituting $f(t)$ into \eqref{inequality}, we have
\begin{align}
 (1-\beta-\alpha)^2 
 \frac{1}{(1-\alpha)\alpha}
\leq
\chi^2(P||Q)
.  \end{align}
Substituting $f(t)$ into \eqref{inequality_symmetric}, we have
\begin{align}
(1-\beta-\alpha)^2  \frac{1}{(1-\beta)\beta}
\leq
\chi^2(Q||P)
.  \end{align}
\end{example}

\subsection{Proof of Proposition \ref{thm_bayes}}

\begin{proof}
$(i) \rightarrow (ii)$: Denote  by $(\alpha^*, \beta^*)$ the tangent point of the supporting line. We have $c = \pi_p \alpha^* + \pi_q \beta^*$. Since BER is the lowest possible value and $c$ is attainable,  we have BER $\leq c$ under class probabilities $(\pi_p,\pi_q)$. Additionally, there is no smaller $c$ that can be attained by $(\alpha,\beta)$ in the Neyman-Pearson region, because $c = \pi_p \alpha + \pi_q \beta$ is the supporting line, meaning that $c \leq \pi_p \alpha + \pi_q \beta$ for all $(\alpha,\beta)$ in the Neyman-Pearson region. Thus, BER $=c$.

$(ii) \rightarrow (i)$: By definition of BER, BER $=c \leq \pi_p \alpha + (1-\pi_p)\beta$ for all $(\alpha,\beta)$ in the Neyman-Pearson region, and the equality is attained by $(\alpha^*, \beta^*)$. Thus, $c \leq \pi_p \alpha + (1-\pi_p)\beta$ is a lower bound of the Neyman-Pearson region and intersects with the region at $(\alpha^*, \beta^*)$. 

$(iii) \leftrightarrow (iv)$ and $(v) \leftrightarrow (vi)$ can be proved analogously. In $(iii) \leftrightarrow (iv)$, since the Neyman-Pearson region is convex, a supporting line to a lower bound of the region is a lower bound for  the entire region. Thus BER $\geq c$ is valid.
\end{proof}

\subsection{Proof of Theorem \ref{thm_refined_upper}}

\begin{proof}
The line of ignorance must be contained in the Neyman-Pearson region, and above the Neyman-Pearson boundary. Any point above the convex envelope is on a segment between a point on the upper bound and $(0,1)$ or $(1,0)$, meaning that it can be attained by a randomized test. Thus, a sharper upper bound can be obtained by the convex envelope.
\end{proof}

\subsection{Proof of Proposition \ref{thm_upper}}

\begin{proof}
Let
\begin{equation}
F(\alpha,\beta,s) = (2-s)\alpha + s\beta 
-
 s^q (2-s)^{1-q} \rho_q
.  \end{equation}
The following two equations hold for the envelope of the family of straight lines \citep[Chapter III]{courant2011differential}
\begin{align}
F(\alpha,\beta,s) 
&=
   (2-s)\alpha + s\beta 
-
 s^q (2-s)^{1-q} \rho_q 
 = 0 \label{eq_f0}
 \\
  \frac{\partial F}{\partial s} (\alpha,\beta,s) 
  &=  
  -\alpha + \beta - s^{q-1}(2-s)^{-q}  (2q-s) \rho_q
  =
  0 \label{eq_pf0}
.  \end{align}
Substitute $\beta$ in \eqref{eq_f0}
 using \eqref{eq_pf0},
 \begin{align}\label{eq_alpha_impl}
 (2-s)\alpha + s \left( \alpha + s^{q-1}(2-s)^{-q}  (2q-s) \rho_q
 \right)
-
  s^q (2-s)^{1-q} \rho_q = 0
.  \end{align}
We obtain an expression of $\alpha$ from \eqref{eq_alpha_impl}
\begin{align}
2\alpha = 
s^q (2-s)^{1-q} \rho_q - s^{q}(2-s)^{-q}  (2q-s) \rho_q
.  \end{align}
That is,
\begin{align}\label{eq_alpha}
\alpha = s^q (2-s)^{-q}\rho_q (1 - q)
.  \end{align}
Substitute $\alpha$ in \eqref{eq_pf0} using \eqref{eq_alpha},
\begin{align}
 \beta = &  s^q (2-s)^{-q}\rho_q (1 - q)  +  s^{q-1}(2-s)^{-q}  (2q-s) \rho_q
 \\
 =&
  s^{q-1}\rho_q (2-s)^{-q} (s (1 - q)  +   (2q-s))
  \\
  =&
    s^{q-1}\rho_q (2-s)^{-q} (2q-qs)
    \\
  =&
  s^{q-1}(2-s)^{1-q} \rho_q q
.  \end{align}
To derive the relation between $\alpha$ and $\beta$, we calculate
\begin{align}
\frac{\beta^q}{\alpha^{q-1}}
=
\frac{s^{q(q-1)}(2-s)^{q(1-q)} \rho_q^q q^q}{s^{q(q-1)} (2-s)^{-q(q-1)}\rho_q^{q-1} (1 - q)^{q-1}}
=
 \rho_q q^q (1 - q)^{1-q}
.  \end{align}
That is
\begin{align}
\beta = \left( \rho_q q^q (1 - q)^{1-q} \right)^{\frac{1}{q}} \alpha^{\frac{q-1}{q}}
.  \end{align}
\end{proof}

\subsection{Proof of Theorem \ref{thm_refined_upper_closed}}

\begin{proof} 
The upper bound is already convex. We consider the envelope formed with the line of ignorance. The tangent lines of the upper bound in Proposition \ref{thm_upper} that go through $(0,1)$ and $(1,0)$ are $\beta =  -\rho_q^{-\frac{1}{q}} \alpha  + 1$ and $\beta = -\rho_q^\frac{1}{q} \alpha + \rho^\frac{1}{q}$, respectively. Any point on the segment between $(0,1)$ and the tangent point can be realized by a randomized test. The case of $(1,0)$ is similar.  
By Theorem \ref{thm_refined_upper}, we obtain the refined bound.
\end{proof}

\subsection{Proof of Theorem \ref{thm_approx}}

\begin{proof} 
Applying the Neyman-Pearson Lemma, the Neyman-Pearson bound of two categorical distributions is piecewise linear, and the $j$-th change point $(\alpha_j, \beta_j)$  is realized by the set of items that have the the largest $j$ density ratios.  We need to verify that the slopes are $\{\frac{p_i}{q_i}\}_{i=1}^n$. 
 We have $\alpha_j = \sum_{i=1}^j  p_i$, $\beta_j = 1- \sum_{i=1}^j  q_i$. Thus the slope is $k_i =- \frac{\beta_i - \beta_{i-1}}{\alpha_i - \alpha_{i-1}} = \frac{q_j}{p_j}$.
\end{proof}

\subsection{Proof of Lemma \ref{thm_convex_cdf}}

\begin{proof}
 By Radon-Nikodym, there is an increasing \( f \) such that
\[
Q(E) = \int_E f dP
\]
for all measurable sets. Hence for any measurable sets \( E \) and \( E' \) which have \( P(E) = P(E') \) consider the sets \( S = E \setminus (E \cap E') \) and \( S' = E' \setminus (E \cap E') \), which have equal \( P \) measure, and also have \( \text{ess sup } S \leq \text{ess inf } S' \), then
\[
Q(S) \leq (\text{ess sup } S) P(S) \leq (\text{ess inf } S') P(S') = (\text{ess inf } S') P(S') \leq Q(S'),
\]
hence
\[
Q(E) = Q(S) + Q(E \cap E') = Q(S') + Q(E \cap E') \leq Q(E').
\]
The ‘leftmost’ subset of \( (0,1) \) of measure \( P(E) \) is \( (0, \mu(E)) \), hence
\[
\alpha((0,\mu(E))) \leq \alpha(E).
\]
\end{proof}

\subsection{Proof of Theorem \ref{thm_exact_realization}}

\begin{proof}
$B^{-1}(1-x)$ is a convex function of $x$. By Lemma \ref{thm_convex_cdf}, the sets that characterize the Neyman-Pearson boundary are intervals $[0, t]$, $t\in [0,1]$. We have $\alpha = \int_0^t dQ = B^{-1}(1-t)$, $\beta = \int_t^1 dP  = 1-t$. Thus, the boundary $\beta = B(\alpha)$ is realized. 
\end{proof}

\subsection{Proof of Proposition \ref{thm_roc_realization}}

\begin{proof}
Consider any point $(t, g(t))$ on the ROC curve. There is a point $(t, f(t))$ on the Neyman-Pearson boundary that has the same horizontal coordinate $t$. This point may not be unique, since the Neyman-Pearson boundary can be a vertical segment at $\alpha=0$. In such case, we take the lowest point. There is another point $(t, 1-t)$ on the line of ignorance, realized by a randomized test.

 Then, the point $(t, g(t))$ can be realized by randomizing the two tests that realize $(t, f(t))$ and $(t, 1-t)$. Specifically, it applies the test realizing $(t, f(t))$  with probability $\frac{1-t-f(t)}{g(t) - f(t)}$, and applies the test realizing $(t, 1-t)$ with probability $\frac{1-t-f(t)}{1- g(t)}$. We have a one-parameter family of tests, where the parameter $t$ determines the tests to be randomized and the probabilities. 
\end{proof}

\end{document}